\documentclass{amsart}
\let\<=\langle
\let\>=\rangle
\let\noi=\noindent
\let\sse=\subseteq
\let\vphi=\varphi

\let\limply=\Longrightarrow

\def\0{\{0\}}

\def\asc{{\rm asc\kern1pt}}
\def\dsc{{\rm dsc\kern1pt}}

\def\esup{{\kern1pt{\rm ess}\kern.5pt\sup\kern1.5pt}}

\def\conv{{\;\longrightarrow\;}}
\def\wconv{{{\buildrel_{\scriptstyle w}\over\conv}}}

\font\nrm=cmr9
\font\esc=cmcsc10 scaled 800

\def\sslash{\hbox{{\fiverm /}}}
\def\notwconv{{{\wconv\kern-13pt\sslash}\kern9pt}}

\def\A{{\mathcal A}}
\def\B{{\mathcal B}}

\def\H{{\mathcal H}}

\def\Le{{\mathcal L}}
\def\M{{\mathcal M}}
\def\N{{\mathcal N}}

\def\R{{\mathcal R}}

\def\X{{\mathcal X}}

\def\BH{{\B[\H]}}

\def\CC{{\mathbb C\kern.5pt}}
\def\FF{{\mathbb F\kern.5pt}}
\def\QQ{{\mathbb Q\kern.5pt}}
\def\NN{{\mathbb N\kern.5pt}}
\def\RR{{\mathbb R\kern.5pt}}
\def\TT{{\mathbb T\kern.5pt}}
\def\ZZ{{\mathbb Z\kern.5pt}}

\newsymbol\varnothing 203F

\newsymbol\blacksquare 1004

\def\smallmatrix#1{\null\,\vcenter{
                   \baselineskip=8pt\mathsurround=0pt\ialign{
                   \hfil ${\scriptstyle##}$
                   \hfil &&
                   \hfil ${\scriptstyle##}$
                   \hfil \crcr
                   \mathstrut \crcr
                   \noalign{\kern-\baselineskip}#1 \crcr
                   \mathstrut \crcr
                   \noalign{\kern-\baselineskip} \crcr }}\!}

\theoremstyle{definition}

\begin{document}

\vglue-80pt\noi
\hfill{\it }\phantom{{\bf XX} (20XX) xx--xx (XXXX) to appear}

\vglue25pt
\title[Posinormality and the Root Problem]
      {Posinormality and the Root Problem}
\author{C.S. Kubrusly}
\address{Catholic University of Rio de Janeiro, Rio de Janeiro, Brasil}
\email{carlos@ele.puc-rio.br}
\author{H.M. Stankovi\'c}
\address{Faculty of Electronic Engineering, University of Ni\v s, Ni\v s,
Serbia}
\email{hranislav.stankovic@elfak.ni.ac.rs}
\renewcommand{\keywordsname}{Keywords}
\date{January 11, 2026}

\vskip-5pt
\begin{abstract}
The paper extends three results regarding the $n$th root problem~by embedding
classes of Hilbert-space operators into the class of posinormal operators.\
For instance, it is shown that (i) for coposinormal operators, if $T\kern-1pt$
is paranormal and $T^n\!$ is quasinormal, then $T\kern-1pt$ is normal, and
(ii) for posinormal operators, if $T\kern-1pt$ is $k$-quasiparanormal and
$T^n\!$ is normal, then $T\kern-1pt$ is normal.\ Moreover, (iii) it is also
shown that the latter result is not conditioned~to~the separability of the
underlying Hilbert space, even if $T\kern-1pt$ is not posinormal, where, in
such a case, $T\kern-1pt$ is the direct sum of a normal operator with a
nilpotent~one.\
\end{abstract}

\maketitle

\vskip-15pt\noi
\section{Introduction}
The so-called $n$th root problem for powers of normal operators on Hilbert
space asks for a description of classes of Hilbert-space operators
$T\kern-1pt$ for which, if $T^n\!$ is~normal, then $T\kern-1pt$ is normal as
well.\ A systematic investigation of the problem dates back to the 1960s and
continues these days.\ Significant examples from the literature will be
discussed here as long as those classes of Hilbert-space operators are
defined.\

\vskip6pt
We bring posinormal operators to the $n$th root problem.\ Classes
of~opera\-tors that are usually investigated towards the problem will be
intercepted with the class of posinormal operators, thus extending some
current results on the $n$th root problem.\

\vskip6pt
Posinormal is a very large class of Hilbert-space operators that includes all
invertible operators and also the hyponormal operators, but not the paranormal
operators.\ $\kern-1pt$The various classes of operators associated with the
class of posinormal operators, for instance, quasiposinormal operators, are
summarised in Section 4.\

\vskip6pt
The central results proved here, related to the above-mentioned classes,
are~$\kern-.5pt$The\-orems 5.4, 5.6, 5.8, and 6.5.\ $\kern-1pt$In particular,
the following corollaries are \hbox{established}.\
\vskip4pt
\begin{description}
\item{$\kern-2pt\circ\kern2pt$}
Suppose $T^*\kern-1pt$ is posinormal.\ If $T\kern-1pt$ is paranormal and
$T^n\kern-1pt$ is quasinormal for some positive integer $n$, then $T\kern-1pt$
is normal (Corollary 5.5).\
\vskip4pt
\item{$\kern-2pt\circ\kern2pt$}
If\/ $T^*\!$ is posinormal and\/ $T^n\kern-1pt$ is hyponormal for some
positive integer $n$, then ${\N(T^n)=\N(T^{*n})}$ for every positive integer
$n$ (Corollary 5.7).\
\vskip4pt
\item{$\kern-2pt\circ\kern2pt$}
Suppose $T\kern-1pt$ is posinormal.\ If\/ $T\kern-1pt$ is $k$-quasiparanormal
and $T^n\kern-1pt$ is normal for~some positive integer $n$, then $T\kern-1pt$
is normal (Corollary 6.6).\
\end{description}
\vskip-3pt

\vskip6pt
Moreover, it is shown that a separability assumption recently required in
\cite[Theorem 3.2]{SK} can be dismissed, yielding the following sharper
result.
\begin{description}
\vskip4pt
\item{$\kern-2pt\circ\kern2pt$}
If $T\kern-1pt$ is $k$-quasiparanormal, and $T^n\kern-1pt$ is normal, then
${T\kern-1pt=N\oplus L}$, where $N$ is normal and $L$ is nilpotent of index at
most $\min\{n,k\kern-1pt+\kern-1pt1\}$ (Theorem 6.3).\
\end{description}
The special case of ${n=2}$ is revisited in Theorem 6.4, and the case where
a $k$-quasi\-paranormal operator has dense range (or is invertible, as a
particular case of being posinormal and coposinormal) is considered in
Theorem 6.7.\

\vskip6pt
The paper is organised into 5 more sections.\ Basic terminology and
notation~are posed in Section 2.\ $\kern-1pt$The classical chain of classes
of operators, from normal to normaloid operators, is considered in Section 3.\
Section 4 gives a short survey on posinormal operators.\ Posinormality and the
$n$th root problem are investigated in Section 5.\ The case of
$k$-quasiparanormal operators is discussed in Section 6.\

\section{Basic Terminology and Notation}

Let $\H$ be a Hilbert space.\ Except in one example, $\H$ is
infinite-dimensional all over the paper.\ If $\M$ is a linear manifold of
$\H$, then $\M^-\!$ denotes the closure of $\M$ in $\H$, and $\M^\perp\!$
is the orthogonal complement of $\M$ in $\H.$ A subspace of $\H$ is a closed
linear manifold of $\H.$ By an operator we mean a bounded linear
transformation of $\H$ into itself.\ Let $\BH$ denote the $C^*\!$-algebra of
all operators on $\H.$ If ${T\kern-1pt\in\BH}$, then ${T^*\!\in\BH}$ stands
for the adjoint of $T\kern-1pt.$ A nonzero operator $T\kern-1pt$ is nilpotent
of index $j$ if ${T^j\kern-1pt=O}$ for some integer ${j>1}$ and ${T^i\!\ne O}$
for every positive integer ${i<j}$, where $O$ denotes the null operator.\ A
nilpotent operator without a specified index is of index 2.\ $\kern-1pt$The
identity operator is denoted by $I$, and $\oplus$ stands for orthogonal direct
sum, either of subspaces or of operators.\ Let ${\N(T)=T^{-1}(\0)}$ and
${\R(T)=T(\H)}$ be the kernel and range of $T\kern-1pt$, respectively; $\N(T)$
is a subspace of $\H$ and $\R(T)$ is a linear manifold of $\H.$ An operator
${T\kern-1pt\in\BH}$ is invertible if it is injective (${\N(T)=\0}$) and
surjective (${\R(T)=\H}).$ Its inverse $T^{-1}\!$ lies in $\BH$ by the Open
Mapping Theorem, and ${T\kern-1pt\in\BH}$ is quasiinvertible if it is
injective and has a dense range ($\R(T)^-\!=\H)$; equivalently, if
$T\kern-1pt$ and $T^*\!$ are~injective, so that $T\kern-1pt$ and $T^*\!$ are
quasiinvertible together.\

\section{The Chain from Normal to Normaloid}

A Hilbert-space operator $T\kern-1pt$ is normal if it commutes with its
adjoint (i.e., $T^*T\kern-1pt=T\kern1pt T^*\kern-1pt).$ In a Hilbert-space
setting, an operator $\kern-.5pt T\kern-.5pt$ is an isometry if and only if\/
${T^*T\kern-1pt=I}$ and $T\kern-1pt$ is a coisometry if $T^*\!$ is an
isometry.\ $\kern-1pt$A unitary operator is an invertible~isometry, that is, a
surjective isometry; equivalently, a normal isometry, which means an isometry
and a coisometry (i.e., ${T^*T\kern-1pt=T\kern1pt T^*\!=I}).$ $\kern-1pt$An
operator $T\kern-1pt$ is self-adjoint if it coincides with its adjoint (i.e.,
${T^*\!=T}).$ $\kern-1pt$Therefore,
$$
\hbox
{\esc Unitary $\scriptstyle\cup$ Self-adjoint $\scriptstyle\subset$ Normal}.
$$
It is quasinormal if it commutes with ${T^*T\kern-1pt}$ (i.e., if
${(T^*T\kern-1pt-T\kern1pt T^*)\kern1pt T\!=O}$), subnormal if it is the
restriction of a normal operator to an invariant subspace (i.e., if it has a
normal extension), hyponormal if ${T\kern1pt T^*\!\le T^*T}$ (i.e.,
${T^*T\kern-1pt-T\kern1pt T^*\!\ge O}$), paranormal if
${\|Tx\|^2\le\|T^2x\|\,\|x\|}$ for \hbox{every} ${x\in\H}$, and normaloid if
${\|T^n\|=\|T\|^n}\!$ for all positive integers $n$ (equivalently,
${r(T)=\|T\|}$, where $r(T)$ stands for the spectral radius of $T).$
$\kern-1pt$These classes are related by proper inclusions, leading to the
classical chain~from normal to normaloid operators,
$$
\kern-4pt\hbox{\esc Normal $\scriptstyle\subset$ Quasinormal
$\scriptstyle\subset$ Subnormal $\scriptstyle\subset$ Hyponormal
$\scriptstyle\subset$ Paranormal $\scriptstyle\subset$ Normaloid}.
                                                               \eqno{(\star)}
$$
Isometries are quasinormal.\ If the adjoint of an operator is quasinormal,
subnormal, hyponormal, or paranormal, then the operator is called
coquasinormal, cosubnormal, cohyponormal, or coparanormal, respectively.\

\section{A Summary of Posinormal Operators}

Definition: an operator $T\kern-1pt$ acting on a Hilbert space is
$$
\hbox{{\it posinormal}$\;$ if $\;{\R(T)\sse\R(T^*)}$,}
$$
\vskip-2pt\noi
$$
\hbox{{\it quasiposinormal}$\;$ if $\;{\R(T)^-\!\sse\R(T^*)^-}\!.$}
$$
\vskip-2pt

\vskip6pt
It is {\it coposinormal}\/ or {\it coquasiposinormal}\/ if its adjoint is
posinormal or quasiposinormal, respectively.\ $\kern-1pt$Thus, $T\kern-1pt$ is
$$
\hbox{coposinormal$\;$ if and only if $\;{\R(T^*)\sse\R(T)}$,}
$$
\vskip-2pt\noi
$$
\hbox{coquasiposinormal$\;$ if and only if $\;{\R(T^*)^-\!\sse\R(T)^-}\!.$}
$$
\vskip2pt\noi
Consequently, $T\kern-1pt$ is
$$
\hbox{posinormal and coposinormal$\;$ if and only if $\;\R(T)=\R(T^*),$}
$$
\vskip-2pt\noi
$$
\hbox{quasiposinormal and coquasiposinormal$\;$ if and only if
$\;\R(T)^-=\R(T^*)^-\!$}.\
$$

\vskip6pt
Concerning the above definitions, recall that $\,{\0^\perp\!=\H}$,
$\,{\Le\sse\M\limply\M^\perp\!\sse\Le^\perp\!}$,
$$
\hbox{${\R(T^*)^\perp\!=\N(T)}$,
\quad
${T\kern-1pt=T^{**}}\!$,
\quad
${\M^-\!=\M^{\perp\perp}\!}$,
\quad and \quad
${\N(T)=\N(T)^-}\!$}
$$
for every operator $T\kern-1pt$ on $\H$ and all linear manifolds $\Le$ and
$\M$ of $\H$, so that
$$
\R(T)\sse\R(T^*)
\;\;\limply\;\,
\R(T)^-\!\sse\R(T^*)^-
\!\iff
\N(T)\sse\N(T^*),                                               \eqno{(\dag)}
$$
and these expressions hold if inclusions are replaced with equalities.\
$\kern1pt$Therefore,
\begin{description}
\item{$\kern-6pt$\rm(i)$\kern2.5pt$}
posinormal operators are quasiposinormal,
\vskip4pt
\item{$\kern-8pt$\rm(ii)$\kern1.5pt$}
$T\kern-1pt$ is quasiposinormal if and only if
${\N(T)\sse\N(T^*)}$,
\vskip4pt
\item{$\kern-9pt$\rm(iii)}
$\kern-.5pt T\kern-1pt$ is coquasiposinormal if and only if
${\N(T^*)\sse\N(T)}$,
\vskip4pt
\item{$\kern-9pt$\rm(iv)}
$T\kern-1pt$ is quasiposinormal and coquasiposinormal if and only
if $\kern-1pt{\N(T)\kern-1.5pt=\kern-1.5pt\N(T^*)}$,
\vskip4pt
\item{$\kern-8pt$\rm(v)$\kern1.5pt$}
the classes of posinormal and quasiposinormal coincide if $\R(T)$ is closed.\
\end{description}

\vskip6pt
Posinormal is a very large class of Hilbert-space operators.\ For instance,
$$
\hbox{injective operators are quasiposinormal,}
$$
\vskip-2pt\noi
$$
\hbox{surjective operators are coposinormal,}
$$
\vskip-0pt\noi
$$
\hbox{invertible operators are posinormal and coposinormal,}
$$
\vskip-0pt\noi
$$
\hbox{quasiinvertible operators are quasiposinormal and coquasiposinormal.}
$$
\vskip-2pt

\vskip6pt
The class of posinormal operators was proposed and investigated by Rhaly
\cite{Rha} in 1994.\ It was observed that every hyponormal operator is
posinormal \cite[Corollary 2.1]{Rha}.\ In fact, posinormal operators include
the hyponormal operators but are not included in the class of normaloid
operators, and paranormal operators are not necessarily posinormal (these
noninclusions were shown in
\cite[Section 5, p$.\,$137]{KD}$\kern.5pt).$~\hbox{Indeed},
$$
\hbox{\esc Hyponormal $\scriptstyle\subset$ Posinormal},
$$
\vskip-2pt\noi
$$
\hbox{\esc Paranormal $\scriptstyle\not\subset$ Posinormal
\qquad {\rm and} \qquad
Posinormal $\scriptstyle\not\subset$ Normaloid}.
$$
So if posinormal operators are inserted in chain $(\star)$, then it goes as
far as this:
$$
\kern-4pt\hbox{\esc Normal $\!\scriptstyle\subset\!$ Quasinormal
$\!\scriptstyle\subset\!$ Subnormal $\!\scriptstyle\subset\!$ Hyponormal
$\!\scriptstyle\subset\!$ Posinormal $\!\scriptstyle\subset\!$
Quasiposinormal}.\!                                               \eqno{(\S)}
$$
\vskip-2pt

\vskip6pt
There are many equivalent definitions of posinormal operators; see, for
instance, in chronological order, \cite[Theorem 2.1]{Rha},
\cite[Theorem B]{Ito}, \cite[Theorem 1]{JKKP}, \cite[Proposition 1]{KD},
\cite[Definition 1]{KVZ}.\ For a comprehensive exposition on posinormal
operators, see, e.g., \cite{Rha} and \cite{KD}.\ The next proposition pulls
together a collection of basic properties of posinormal operators that will
be required throughout~the~\hbox{paper}.\

\vskip6pt\noi
{\bf Proposition 4.1.}\
Let $T\kern-1pt$ be an operator acting on a Hilbert space.\
\begin{description}
\item{$\kern-12pt$\rm(a)$\kern3pt$}
If $T\kern-1pt$ is posinormal, then $\N(T)$ reduces $T\;$
\cite[Corollary 2.3]{Rha}.
\vskip4pt
\item{$\kern-12pt$\rm(b)$\kern2.5pt$}
$T\kern-1pt$ is quasiposinormal if and only if $\N(T)$ reduces $T\;$
(proof of \cite[Lemma 1]{KD}).\
\vskip4pt
\item{$\kern-11pt$\rm(c)$\kern2pt$}
Restrictions of posinormal operators to invariant subspaces are posinormal$\;$
\hbox{\cite[Proposition 4]{JKKP}}.\
\vskip4pt
\item{$\kern-12pt$\rm(d)$\kern3.5pt$}
If $T\kern-1pt$ has no nontrivial invariant subspace, then it is
quasiposinormal \hbox{\cite[$\kern-1pt$p$.\kern1pt$139]{KD}}.\
\vskip4pt
\item{$\kern-11pt$\rm(e)$\kern2.5pt$}
If $T\kern-1pt$ is quasiposinormal, then ${\N(T^2)=\N(T)}$ --- see
\cite[Proposition 3]{JKKP} (also, \cite[Remark 2]{KD}$\kern.5pt$) for
posinormal operators; the same proof holds for quasiposinormal operators ---
so if $T$ is quasiposinormal, then ${\N(T^n)=\N(T)}$~\hbox{for all $n.$}
\vskip4pt
\item{$\kern-12pt$\rm(f)$\kern3pt$}
$T\kern-1pt$ being posinormal does not imply that $T^2$ is posinormal$\;$
\cite[Example 1]{KVZ}.\
\vskip4pt
\item{$\kern-12pt$\rm(g)$\kern3pt$}
If $T\kern-1pt$ is quasiposinormal, then $T^n$ is quasiposinormal for every
positive \hbox{integer} $n\;$ \cite[Theorem 5]{KVZ}.\
\vskip4pt
\item{$\kern-12pt$\rm(h)$\kern1pt$}
If $T\kern-1pt$ is posinormal and $\R(T)$ is closed, then $T^n\kern-1pt$ is
posinormal and $\R(T^n)$ is closed for every positive integer $n\;$
\cite[Corollary 14]{JV} (also, \cite[Main~Theorem]{BKT}$\kern.5pt$).\
\vskip4pt
\item{$\kern-11pt$\rm(i)$\kern3pt$}
If $T\kern-1pt$ has a closed range and is hyponormal, then $T^n$ has a closed
range for every positive integer $n$ \cite[Corollary 2.3]{BKT}.\
\end{description}

\vskip6pt
As closed-range posinormal operators are quasiposinormal, (h) is a nontrivial
special case of (g) and (i) is a consequence of (h).\ Item (e) also holds for
paranormal operators.\ The product of closed-range posinormal operators was
investigated in \cite{BKLT}.\

\vskip6pt\noi
{\bf Remark 4.2.}\
An operator and its adjoint are normal or normaloid together$.$~Also,
\vskip6pt\noi
{\narrower
if an operator is hyponormal and cohyponormal, then it is normal,
\vskip0pt}
\vskip4pt\noi
and so if it is subnormal and cosubnormal, or quasinormal and coquasinormal,
then it is normal.\ The first attempts to get a result along this line for
paranormal operators involved posinormality assumptions.\ It was shown in
\cite[Theorem 5]{And} that
\vskip6pt\noi
{\narrower
{if $T\kern-1pt$ is paranormal and coparanormal and
${\N(T)=\N(T^*\kern-.5pt)}$, then $T\kern-1pt$ is normal.}
\vskip0pt}
\vskip6pt\noi
This was followed by \cite[Theorem 8.3.47]{Ist} with a similar proof, although
with a slightly more restrictive condition (cf.\ $(\dag)\kern.5pt$), namely,
\vskip6pt\noi
{\narrower
{if $T\kern-1pt$ is paranormal and coparanormal and
${\R(T)=\R(T^*\kern-.5pt)}$, then $T\kern-1pt$ is normal.}
\vskip0pt}
\vskip6pt\noi
Note that these results can be rewritten, respectively, as follows:
\vskip6pt\noi
{\narrower
{if an operator is paranormal and coparanormal, and if it is quasiposinormal
and coquasiposinormal, then it is normal,} \qquad and
\vskip0pt}
\vskip4pt\noi
{\narrower
{if an operator is paranormal and coparanormal, and if it is posinormal
and coposinormal, then it is normal.}
\vskip0pt}
\vskip6pt\noi
However, the above posinormality conditions can be dismissed:
\vskip6pt\noi
{\narrower
{if an operator is paranormal and coparanormal, then it is normal,}
\vskip0pt}
\vskip6pt\noi
as was shown later in \cite[Corollary 3]{YY}.\

\section{Posinormality and the $n$th Root Problem}

Throughout this section, $\H$ stands for a Hilbert space, and operators act on
Hilbert spaces.\ A first consequence of the above definitions reads as follows.

\vskip6pt\noi
{\bf Proposition 5.1.}\
{\it If\/ $T\kern-1pt$ is quasiposinormal and coquasiposinormal, then
$$
T\kern-1pt=O\oplus S
\quad\;\hbox{on}\;\quad
\H=\M\oplus\M^\perp\!,
$$
where\/ ${\M=\N(T)}$, $\,{\M^\perp\!=\R(T)^-}\!$, and\/
${S=T|_{\R(T)^-}}\!$ is quasiinvertible}\/.\

\proof
Let $T\kern-1pt$ be quasiposinormal and coquasiposinormal (in particular,
posinormal and coposinormal) so that ${\R(T)^-\!=\R(T^*)^-}\!$,
equivalently, ${\N(T)=\N(T^*)}.$ Since ${\N(T)=\R(T^*)^\perp}\!$ for every
operator, ${\N(T)=\R(T^*)^\perp\!=\R(T)^\perp}\!.$ By
Proposition 4.1(b), ${\N(T)=\R(T)^\perp}\!$ reduces $T\kern-1pt$, and so
${\R(T)^-\!=\R(T)^{\perp\perp}}\!$ reduces $T\kern-1pt.$~$\kern-1pt$Thus
$$
T\kern-1pt=O\oplus T|_{\R(T)^-}
\quad\hbox{on}\;\quad
\H=\N(T)\oplus\R(T)^-\!.
$$
Since $\R(T)^-\!$ reduces $T\kern-1pt$ and is orthogonal to $\N(T)$,
$\,{T|_{\R(T)^-}\!:\R(T)^-\!\!\to\R(T)}\sse\R(T)^-\!$ has a dense range,
${\R(T|_{\R(T)^-\kern-1pt})^-\!=\R(T)^-}\!$, and is injective,
$\N(T|_{\R(T)^-\kern-1pt})^-\!=\0$, which means that
${T|_{\R(T)^-\!}\in\B[\R(T)^-\!]}$ is quasiinvertible.                   \qed

\vskip6pt
An injective operator $T\kern-1pt$ is quasiposinormal and coquasiposinormal
if and only if ${\N(T)=\N(T^*)=\0}.$ This means that $T\kern-1pt$ is
quasiinvertible.\ In this case, the null part $O$ is absent in the orthogonal
direct sum decomposition of Proposition 5.1.\

\vskip6pt
Proposition 5.1 leads us to the following result.\ {\it If a quasinormal
operator is~co\-quasiposinormal, then it is normal}\/ \cite[Lemma 2.4]{SK}.\
Such a result will play an~important role in the sequel.\ Actually, Lemma 2.4
in \cite{SK} reads as follows:\ {\it for quasi\-normal operators,\/ {\rm(b)}
implies}\/ {\rm(a)}, as a part of the next lemma.\

\vskip6pt\noi
{\bf Lemma 5.2.}\
{\it Take an operator $T\kern-1pt$ and consider the following assertions.
\begin{description}
\item{$\kern-12pt$\rm(a)$\kern0pt$}
$T\kern-1pt$ is normal\/.
\vskip4pt
\item{$\kern-12pt$\rm(b)$\kern1pt$}
$T\kern-1pt$ is coquasiposinormal
$\hfill$
$($i.e.,\/ ${\N(T)\sse\N(T^*)}$, or\/ ${\R(T^*)^-\!\sse\R(T)^-\!})$.
\vskip4pt
\item{$\kern-12pt$\rm(c)$\kern.5pt$}
$\R(T)^-\!$ reduces\/ $T$
$\kern62pt$
$($equivalently,\/ $\N(T^*)$ reduces\/ $T)$.
\end{description}
\vskip1pt\noi
For an arbitrary operator\/ $T\kern-1pt$,\/
${{\rm(a)}\Rightarrow{\rm(b)}\Rightarrow{\rm(c)}}.$ If\/ $T\kern-1pt$ is
quasinormal, then\/ ${{\rm(c)}\Rightarrow{\rm(a)}}.$
$\kern-1pt$\hbox{Therefore}, if\/ $\kern-.5pt T\kern-1.5pt$ is quasinormal,
then the above assertions are pairwise equivalent}\/.\

\proof
${{\rm(a)}\Rightarrow{\rm(b)}}.$
Since the adjoint of a normal operator is normal, (a) implies (b) according to
the chain in $(\S).$
\vskip6pt\noi
${{\rm(b)}\Rightarrow{\rm(c)}}.$
If $T\kern-1pt$ is coquasiposinormal, then $\N(T^*)$ reduces $T^*\!$ by
Proposition 4.1(b) --- or Proposition 5.1 --- and hence it reduces
$T\kern-1pt.$ So (b) implies (c).\
\vskip6pt\noi
${{\rm(c)}\Rightarrow{\rm(a)}}.$
Let $T\kern-1pt$ be a quasinormal operator on a Hilbert space $\H.$ Suppose
$\R(T)^-\!$~reduces $T\kern-1pt.$ Since
${(\R(T)^-)^\perp\!=\R(T)^\perp\!}$ $=\N(T^*)$ for every operator on $\H$,
$$
T\kern-1pt=R\oplus S
\quad\;\hbox{and}\;\quad
T^*\kern-1pt=R^*\!\oplus S^*\kern-1pt
\quad\;\hbox{on}\;\quad
\H=\N(T^*)\oplus\R(T)^-\!.
$$
If $T\kern-1pt$ is quasinormal (i.e.,
${(T^*T\kern-1pt-T\kern1pt T^*)\kern1pt T\kern-1pt=O}$), then
${S\kern-.5pt
=\kern-.5pt T|_{\R(T)^-}\kern-2.5pt:\R(T)^-\kern-2.5pt\to\R(T)}\sse\R(T)^-\!$
is quasinormal:\ ${(S^*S-S\kern1pt S^*)S=O}.$ Thus
${(S^*S-S\kern1pt S^*)\kern1pt u=0}$ for every $u$ in $\R(S)$, and so
${(S^*S-S\kern1pt S^*)\kern1pt u=0}$ for every $u$ in $\R(S)^-\!$ (extension
by continuity)$.$~But
${\R(S)^-\kern-2.5pt=\R(T|_{\R(T)^-\kern-1pt})^-\kern-2.5pt
=\R(T)^-\kern-2.5pt}.$
So $S$ on $\R(T)^-\!$ is normal.\ Also,
${R^*\kern-2.5pt=\kern-.5pt T^*|_{\N(T^*)}\!=O}$, and hence ${R=O}.$ Thus
${T\kern-1pt=O\oplus S}$ is normal.\ Then (c) implies (a).               \qed

\vskip6pt
Since an invertible operator is posinormal and coposinormal, thus
quasiposinor\-mal and coquasiposinormal, and since there are invertible
bilateral subnormal~(nonnormal) weighted shifts (see, e.g.,
\cite[p$.\,$60, Remark 5]{Con}$\kern.5pt$), the equivalence in Lemma 5.2 does
not survive the trip from quasinormal to subnormal (and so
to~\hbox{hyponormal}).\

\vskip6pt
It is a trivial (and widely known) result that an invertible quasinormal
operator is normal.\ Lemma 5.2 leads to immediate generalisations of it:\
\vskip4pt\noi
\centerline{\it a dense-range quasinormal operator is normal\/{\rm;}}
\vskip2pt\noi
\centerline{\it an injective quasinormal operator is normal\/.}

\vskip6pt
The $n$th root problem asks for classes of nonnormal operators $T\kern-1pt$
with the following property:\ if $T^n$ is normal, then $T\kern-1pt$ is
normal.\ The functional calculi for normal operators on complex Hilbert spaces
ensure that a normal operator $N$ has a normal $n$th root.\ For instance, if
${\psi\!:\CC\to\CC}$ is the function assigning to each
${\zeta\kern-1pt\in\CC}$ its principal $n$th root, that is,
${\psi(\zeta)=|\zeta|^{\frac{1}{\scriptstyle n}}
\hbox{exp}\big(i\,\hbox{Arg}(\zeta)\kern1pt\frac{1}{\scriptstyle n}\big)}$,
then
${\psi(\zeta)^n=\zeta}$ for every ${\zeta\kern-1pt\in\CC}$ and $\psi(N)$ is a
normal $n$th root of $N.$ However, there are square roots of normal operators
that are not normal (e.g., every nilpotent operator is a square root of the
null operator, which are never normal unless they are null).\ Also,
there~may~be~several normal square roots of a normal operator (trivial
example:\ $I$ and $-I$ are normal square roots of the identity operator).\ For
more along these lines, see,~e.g.,~\cite{Dug,Ker}.\

\vskip6pt
Originally, part of Lemma 5.2 was posed in \cite[Lemma 2.4]{SK} to give an
alternate and straightforward proof of a simple $n$th root problem considered
in \cite[Lemma 2.3]{SK}, which reads as follows.\ Let $n$ be an arbitrary
positive integer.\
$$
\hbox{\it If\/ $T\kern-1pt$ is quasinormal and\/ $T^n\kern-1pt$ is normal,
then $T\kern-1pt$ is normal}\/.                         \eqno\hbox{\nrm(5.1)}
$$
\vskip2pt\noi
This comes as a particular case, with an elementary proof, of a useful and
well-known result on hyponormal operators, given in~\cite[Theorem 5]{Stampfli}:
$$
\hbox{\it If\/ $T\kern-1pt$ is hyponormal and\/ $T^n\kern-1pt$ is normal,
then\/ $T$ is normal}\/,                                \eqno\hbox{\nrm(5.2)}
$$
which was later extended to paranormal operators in \cite[Theorem 6]{And}:
$$
\hbox{\it If\/ $T\kern-1pt$ is paranormal and\/ $T^n\kern-1pt$ is normal,
then\/ $T$ is normal}\/.                                \eqno\hbox{\nrm(5.3)}
$$
These exhaust the `$T^n$-normal implies $T$-normal' problem among the classes
of~operators exhibited in chain $(\star).$ Of course, there is an uncountable
number of classes of operators that include the class of paranormal operators
and are included in the class of normaloid operators, but we will not deal
with them in this section.\ $\kern-1pt$The problem cannot be extended to
normaloid operators, since normaloid is too large a class.\ For instance, if
a nonzero operator $L$ is a nilpotent contraction (i.e.,~if~${L^2\kern-1pt=O}$
and ${0\kern-1pt<\|L\|\le1}$), then $L$ is not normal (it is not
\hbox{normaloid}),~and~so~${T\kern-1pt=\kern-1pt I\kern-1pt\oplus L}$~is~not
normal.\ However, $T^2\!=$ $\kern-1pt{I\kern-1pt\oplus O}$ is a nonzero
normal, and $T\kern-1pt$ is normaloid with $\|T\|={r(T)=1}$ since
${r(T^2)=1}\,$ (e.g., take
$L\kern-1pt=\kern-1pt\big(\smallmatrix{O & C \cr
                                       O & O \cr}\big)$ with
${0\kern-1pt<\kern-1pt\|C\|\kern-1pt\le1\kern-1pt}.)$ In particular,
$$
T\kern-1pt=\Big(\smallmatrix{1 & 0 & 0 \cr
                             0 & 0 & 1 \cr
                             0 & 0 & 0 \cr}\Big)\;
\hbox
{is a normaloid operator with $T^2\kern-1pt$ normal, but $T$ is not normal}.
$$
In fact, $T\kern-1pt$ is not even paranormal since
${T\!=I\!\oplus\hbox{(Nonzero Nilpotent)}}$ (see also
\cite[Example 2.6]{SK}$\kern.5pt$) and, although with closed range, it is not
posinormal or \hbox{coposinormal}.\

\vskip6pt
An alternative approach regarding the classes of operators in the chain
$(\star)$, beyond the class of quasinormal, is to replace the expression
`$T^n$ normal', as in (5.1)--(5.3) above, with the expression `$T^n$
quasinormal'.\ Recall that if $T\kern-1pt$ is quasinormal, then $T^n$ is
quasinormal for every positive integer $n$ (a consequence of the fact that if
$T\kern-1pt$ is quasinormal, then $T\kern-1pt$ commutes with ${T^{*n}T^n}$ for
every $n$; see, e.g., \cite[Problem 7.2(c)]{HSO}$\kern.5pt$).\ A sort of
converse was recently proved in \cite[Theorem 6.3]{PS1}:
$$
\hbox{\it If\/ $T\kern-1pt$ is subnormal and\/ $T^n\!$ is quasinormal, then\/
$T\kern-1pt$ is quasinormal}.                           \eqno\hbox{\nrm(5.4)}
$$
(For an elementary proof, see \cite[Theorem 2.2 and Corollary 2.3]{Sta},
and for a multivariable version of such a problem, see \cite{Sta2}.)
$\kern-1pt$This has been extended to a class of operators that includes the
hyponormal operators, whose definition runs as follows.\
Set $|T|=(T^*T)^\frac{1}{2}$ and ${D=T^*T\kern-1pt-T\kern1pt T^*}\!.$ Recall
that $T\kern-1pt$~is hyponormal if and only if ${D\ge O}$; equivalently,
${|T|^2\kern-1pt\le|T^*|^2}\!.$ A Hilbert-space~opera\-tor $T\kern-1pt$ is
{\it quasihyponormal}\/ if ${T^*\kern-1pt D\kern1pt T\kern-1pt\ge O}$;
equivalently, ${|T|^4\kern-1pt\le|T^2|^2}\!$, and
{\it semi-quasi\-hyponormal}\/ if ${|T|^2\kern-1pt\le|T^2|}.$
Semi-quasihyponormal operators are also called {\it class}\/~$\A$ operators.\
As is well known, these classes of operators are related by proper inclusions
(see, e.g., \cite[p$.\,$98]{HSO}$\kern.5pt$):
$$
\hbox{\esc Hyponormal$\;\subset\;$Quasihyponormal$
\;\subset\;$Semi-quasihyponormal$\;\subset\;$Paranormal}.      \eqno{(\ddag)}
$$
Like paranormal operators, quasihyponormal (and consequently,
semi-quasihypo\-normal) operators are not posinormal (see
\cite[Example 1]{KD}$\kern.5pt$):
$$
\hbox{\esc Quasihyponormal $\scriptstyle\not\subset$ Posinormal}.
$$
The result in (5.4) above was extended from subnormal to semi-quasihyponormal
operators in \cite[Theorem 4.1]{PS2}:
\kern5pt\noi
$$
\kern-8pt
\hbox{\it If\/ $\kern-.5pt T\kern-1.5pt$ is semi-quasihyponormal and\/
$\kern-.5pt T^n\kern-2pt$ is quasinormal, then\/ $\kern-.5pt T\kern-1.5pt$ is
quasinormal}.                                           \eqno\hbox{\nrm(5.5)}
$$
\vskip0pt
According to $(\ddag)$, the result in (5.5) naturally motivates the following
question.

\vskip6pt\noi
{\bf Question 5.3.}\
{\it Is it true that if\/ $T\kern-1.5pt$ is paranormal and\/ $T^n\kern-2pt$ is
quasinormal, then\/ $T\kern-1.5pt$ is quasinormal}$\,$?

\vskip6pt
However, under the assumption of coquasiposinormality, as in Lemma 5.2, we get
an affirmative answer (partial from one hand and stronger from the other hand).

\vskip6pt\noi
{\bf Theorem 5.4.}\
{\it Suppose\/ $\kern-.5pt T\kern-1.5pt$ is coquasiposinormal.\ If\/
$\kern-.5pt T\kern-1.5pt$ is paranormal and\/ $T^n\!$ is~quasi\-normal for
some positive integer\/ $n$, then\/ $\kern-.5pt T\kern-1.5pt$ is normal}\/.

\proof
Fix an arbitrary positive integer $n.$

\vskip6pt\noi
(a)
According to Proposition 4.1(g), if $T^*\!$ is quasi\-posinormal, then so is
$T^{*n}\!.$ Equivalently, if $T\kern-1pt$ is coquasiposinormal, then
$T^n\kern-1pt$ is coquasiposinormal.\ In this case, if $T^n\kern-1pt$ is
quasinormal, then $T^n\!$ is normal by Lemma~5.2.\
\vskip6pt\noi
(b)
Nevertheless, if $T\kern-1pt$ is paranormal and $T^n\kern-1pt$ is normal, then
$T\kern-1pt$ is normal by (5.3)$.$                                       \qed

\vskip6pt
A particular and simpler form of Theorem 5.4 reads as follows.

\vskip6pt\noi
{\bf Corollary 5.5.}\
{\it If\/ $T\kern-1pt$ is coposinormal, paranormal, and $T^n\!$ is quasinormal
for some positive integer\/ $n$, then\/ $T\kern-1pt$ is normal}\/.

\vskip6pt
As every positive integer power of a quasinormal operator is quasinormal,
we~recover part of Lemma 5.2 as an immediate consequence of Corollary~5.5:\
\vskip4pt
\centerline{\it a coposinormal quasinormal operator is normal\/.}

\vskip6pt
A root problem type for the posinormal family only concludes this section.

\vskip6pt
But first recall that, for an arbitrary operator $A$,
${\N(A^m)\kern-1pt\sse\kern-1pt\N(A^{m+1})}$ for \hbox{every}~integer
${m\kern-1pt\ge\kern-1pt0}$, and if ${\N(A^{m_0})=\N(A^{m_0+1})}$ for some
integer ${m_0\kern-1pt\ge\kern-1pt0}$, then $\N(A^m)=\N(A^{m_0})$ for every
${m\kern-1pt\ge\kern-1pt m_0}.$ Similarly, ${\R(A^{m+1})\sse\R(A^m)}$ for
every ${m\kern-1pt\ge\kern-1pt0}$, and if ${\R(A^{m_0+1})=\R(A^{m_0})}$ for
some integer ${m_0\kern-1pt\ge\kern-1pt0}$, then ${\R(A^m)=\R(A^{m_0})}$ for
\hbox{every} ${m\kern-1pt\ge\kern-1pt m_0}.$ If there exist those integers,
then the least integer $m_0$ for which $\N(A^{m_0})=\N(A^{m_0+1})$ is the
{\it ascent}\/ of $A$, denoted by $\asc(A)$, and the least integer $m_0$ for
which ${\R(A^{m_0+1})=\R(A^{m_0})}$ is the {\it descent}\/ of $A$, denoted by
$\dsc(A).$

\vskip6pt\noi
{\bf Theorem 5.6.}\
{\it If\/ $T\kern-1.5pt$ is coquasiposinormal and $T^n\!$ is quasiposinormal
for some~inte\-ger ${n\kern-1pt\ge\kern-1pt1}$, then\/ $T^n\!$ is
quasiposinormal and coquasiposinormal for every
integer~${n\kern-1pt\ge\kern-1pt1}.$

\proof
As we saw above, ${\N(A^2)=\N(A)}$ implies ${\N(A^m)=\N(A)}$ for all
${m\ge1}.$ Now since $T^*\!$ and $T^n\!$ are quasiposinormal,
${\N(T^{*2})=\N(T^*)}$ and ${\N(T^{2n})=\N(T^n)}$ by Proposition 4.1(e).\
$\kern-1pt$Therefore, for every positive integer $m$,
$$
\N(T^{*n\kern.5pt m})=\N(T^*)
\quad\:\hbox{and}\;\quad
\N(T^{n\kern.5pt m})=\N(T^n).
$$
Moreover, since $T^n\!$ is quasiposinormal, it follows by Proposition 4.1(g)
that $T^{n\kern.5pt m}\!$ is quasiposinormal.\ $\kern-1pt$Therefore, since
$T^*\!$ is also quasiposinormal, we get
$$
\N(T^{n\kern.5pt m})\sse\N(T^{*n\kern.5pt m})
\quad\:\hbox{and}\;\quad
\N(T^*)\sse\N(T).
$$
Hence,
$\;{\N(T)\sse\N(T^n)=\N(T^{n\kern.5pt m})
\sse\N(T^{*n\kern.5pt m})=\N(T^*)\sse\N(T)}\;$
so that
$$
\N(T)=\N(T^*).
$$
Then, besides ${\N(T^{*2})=\N(T^*)}$, as we saw above, we also get
${\N(T^2)=\N(T)}$ by Proposition 4.1(e).\ These imply that, for every
${n\ge1}$,
$$
\N(T^{*n})=\N(T^*)
\quad\;\hbox{and}\;\quad
\N(T^n)=\N(T),
$$
and so ${\N(T^n)=\N(T^{*n})}\,$ (i.e., $T^n\kern-1pt$ is quasiposinormal and
coquasiposinormal).\                                                     \qed

\vskip6pt\noi
{\it Note}\/:\ Observe that {\it if\/ $T\kern-1pt$ is quasiposinormal and\/
$T^n\!$ is coquasiposinormal for some integer\/ ${n\kern-1pt\ge\kern-1pt1}$\/,
then\/ $T^n\!$ is quasiposinormal and coquasiposinormal for every integer\/
${n\kern-1pt\ge\kern-1pt1}$} (since, by Theorem 5.6, $T^{*n}\!$ is
quasiposinormal and coquasiposinormal).\

\vskip6pt
A useful and immediate consequence of Theorem 5.6 reads as follows.

\vskip6pt\noi
{\bf Corollary 5.7.}\
{\it If\/ ${\N(T^*)\sse\N(T)}$ and\/ $T^n\!$ is hyponormal for some positive
integer $n$, then\/ ${\N(T^n)=\N(T^{*n})}$ for every positive integer}\/ $n.$

\vskip6pt
And an extension of Theorem 5.6 is given below (where
$\operatorname{lcm}(k,\ell)$ stands for the least common multiple of $k$ and
$\ell$).\

\vskip6pt\noi
{\bf Theorem 5.8.}\
{\it If\/ $T^k\!$ is quasiposinormal and\/ $T^\ell\!$ is coquasiposinormal
for some pair\/ ${k,\ell}$ of positive integers, then\/ $T^m\!$ is both
quasipo\-sinormal and coquasiposinormal for every\/
${m\kern-1pt\ge\kern-1pt m_0}$, where\/
${m_0\kern-1pt=\kern-1pt\operatorname{lcm}(k,\ell)}.$ Also}\/,
\vskip4pt
\begin{description}
\item{$\kern-6pt$\rm(a)$\kern0pt$}
{\it if\/ ${\operatorname{dsc}(T)\kern-.5pt=\kern-.5pt1}$, then\/ $T\kern-1pt$
is quasiposinormal}\/;
\vskip4pt
\item{$\kern-6pt$\rm(b)$\kern0pt$}
{\it if\/ ${\operatorname{asc}(T)\kern-.5pt=\kern-.5pt1}$, then\/ $T\kern-1pt$
is coquasiposinormal}\/.\
\end{description}

\proof
Let $T^k\!$ be quasiposinormal and\/ $T^\ell\!$ be coquasiposinormal.\ Set
${m_0=\operatorname{lcm}(k,\ell)}.$ By Proposition 4.1(g), it is
straightforward to see that $T^{m_0}\!$ is both quasiposinormal and
coquasiposinormal, that is,
$$
\N(T^{m_0})=\N({T^*}^{m_0}).
$$
Now let $m$ be an arbitrary integer such that ${m\kern-1pt\ge\kern-1pt m_0}.$
Thus, there exists an integer ${p\kern-1pt\ge\kern-1pt1}$ such that
${m\kern-1pt\le\kern-1pt p\,m_0}.$ $\kern-1pt$Then, using Proposition 4.1(e),
the above identity, and the same argument as before (cf.\ proof of
Theorem 5.6), we get
$$
\N(T^m)\sse\N(T^{p\,m_0})=\N(T^{m_0})=\N(T^{*m_0})\sse\N(T^{*m}).
$$
Dually, ${\N(T^{*m})\sse\N(T^m)}.$ Hence, ${\N(T^m)=\N(T^{*m})}$, and so
$T^m\!$ is both quasiposinormal and coquasiposinormal for every ${m\ge m_0}.$
Also, arguing analogously,
\vskip6pt\noi
(a)
If\/ ${\dsc(T)\kern-1pt=\kern-1pt1}$, then ${\N(T^{*m})=\N(T^*)}$ (i.e.,
${\R(T^m)^-\!=\R(T)^-\!}$) for all ${m\kern-1pt\ge\kern-1pt1}.$~As
$T^k\!$ is quasiposinormal, ${\N(T^k)\sse\N(T^{*k})}.$ So
$\,{\N(T)\sse\N(T^k)\sse\N(T^{*k})=\N(T^*)}.$
\vskip6pt\noi
(b)
Again and dually, if\/ ${\asc(T)\kern-1pt=\kern-1pt1}$, then ${\N(T^m)=\N(T)}$
for all ${m\kern-1pt\ge\kern-1pt1}.$ Since $T^\ell\!$ is coquasi\-posinormal,
${\N(T^{*\ell})\sse\N(T^\ell)}.$ So
$\,{\N(T^*)\sse\N(T^{*\ell})\sse\N(T^\ell)=\N(T)}.\!$                    \qed

\section{The $k$-Quasiparanormal Case}

Let $\H$ be a Hilbert space.\ (In fact, the definitions below hold in a normed
space.) Consider the following classes of operators, all of them including
the class of paranormal operators, which have been often discussed and
compared in current literature.\ (See, e.g., \cite[Section 1]{SK} for a
detailed account of these classes$.$) An operator $T\kern-1pt$ on $\H$ is
$k$-{\it paranormal}\/ if ${\|Tx\|^{k+1}\kern-1pt\le\|T^{k+1} x\|\,\|x\|^k}$
for every ${x\in\H}$, for some ${k\ge0}.$ It is {\it quasiparanormal}\/ if it
acts as a paranormal operator on its range; that is, if
${\|T(T x)\|^2\kern-1pt\le\|T^2(T x)\|\,\|T x\|}$ for every ${x\in\H}.$
Similarly, an operator $T\kern-1pt$ on $\H$ is $k$-{\it quasiparanormal}\/ if
it acts as a paranormal operator on the range of $T^k$; that~is,~if
$$
\|T(T^k x)\|^2\le\|T^2(T^k x)\|\,\|T^k x\|
$$
for every ${x\in\H}$, for some positive integer $k.$ In other words, if
$$
\|Ty\|^2\le\|T^2y\|\,\|y\|
$$
for every ${y\in\R(T^k)}$, equivalently (by continuity), for every
${y\in\R(T^k)^-\!}.$

\vskip6pt
It is clear that $\R(T^k)$ is $T\kern-1pt$-invariant, and so is $\R(T^k)^-\!.$
Also, an operator $T$ is surjective if and only if $T^k\!$ is surjective, and
$T$ has dense range if and only if $T^k$ has dense range, for every positive
integer $k$ (i.e., ${\R(T)=\H}$ implies ${\R(T^k)=\H}$ and
${\R(T)^-\kern-3pt=\kern-.5pt\H}$ implies ${\R(T^k)^-\kern-3pt=\kern-.5pt\H}$
for every $k$, which are readily verified by~induction)$.$ These facts and the
above definition of $k$-quasiparanormality are enough to ensure the following
auxiliary results that will be required later in this~section.\

\vskip6pt\noi
{\bf Proposition 6.1.}\
{\it Let\/ $T\kern-1pt$ be an operator on a Hilbert space\/ $\H.$}
\begin{description}
\item{$\kern-9pt$\rm(a)$\kern1pt$}
{\it $T\kern-1pt$ is\/ $\kern-1pt k$-quasiparanormal if and only if\/
$T|_{\R(T^k)^-\!}$ is paranormal on\/ $\R(T^k)^-\!.$}
\vskip4pt
\item{$\kern-9pt$\rm(b)}
{\it If\/ $T$ has a dense range\/ $($in particular, is surjective\/$)$, then
the following asser\-tions are pairwise equivalent.}
\vskip4pt\noi
\begin{description}
\item{$\kern-6pt$\rm(i)$\kern2.5pt$}
{\it $T$ is paranormal.}
\vskip4pt
\item{$\kern-8pt$\rm(ii)$\kern1.5pt$}
{\it $T$ is\/ $k$-quasiparanormal for every\/ $k$.}
\vskip4pt
\item{$\kern-9pt$\rm(iii)}
{\it $T$ is\/ $k$-quasiparanormal for some\/ $k$.}
\end{description}
\end{description}
\goodbreak

\vskip4pt
Also note that
\vskip2pt\noi
\begin{description}
\item{$\bullet\quad$}
$0$-paranormal = $\BH$,
\item{$\bullet\quad$}
$1$-paranormal = paranormal,
\item{$\bullet\quad$}
0-quasiparanormal = paranormal,
\item{$\bullet\quad$}
1-quasiparanormal = quasiparanormal,
\item{$\bullet\quad$}
$k$-paranormal for every $k$ = paranormal,
\item{$\bullet\quad$}
$k$-quasiparanormal $\sse$ ${k+1}$-quasiparanormal,
\end{description}
\vskip2pt\noi
so if $T$ is $k$-quasiparanormal, then $T$ is $m$-quasiparanormal for every
integer ${m\kern-.5pt\ge\kern-1pt k}.$

\vskip6pt
The class of $k$-paranormal operators fits into the chain $(\star)$ of
Section 3, but not the class of $k$-quasiparanormal operators.\ Indeed, as is
well known,
$$
\hbox{\esc Paranormal}\subset k\hbox{\esc-Paranormal}\subset
\hbox{\esc Normaloid}
$$
\vskip-6pt\noi
and
\vskip-3pt\noi
$$
\hbox{\esc Paranormal}\subset k\hbox{\esc-Quasiparanormal}\not\sse
\hbox{\esc Normaloid}.
$$
\vskip1pt\noi
Since paranormal operators are not posinormal, the operators defined in this
section are not posinormal.\

\vskip6pt
The following $n$th root results on $k$-paranormal and $k$-quasiparanormal
operators have been proved in \cite[Theorems 3.1 and 3.2]{SK}.\
$$
\hbox{\it If\/ $T\kern-1pt$ is\/ $k$-paranormal and\/ $T^n\!$ is normal,
then\/ $T\kern-1pt$ is normal}.                         \eqno\hbox{\nrm(6.1)}
$$
$$
\vbox{\hskip-22pt
\vbox{\hskip-12pt
{\narrower\narrower
{\it If\/ $T\kern-1pt$ is a\/ $k$-quasiparanormal operator on a separable
Hilbert space and\/~$T^n\!$ is normal, then\/ ${T\kern-1pt=N\oplus L}$,
where\/ $N$ is normal and\/ $L$ is nilpotent of index at most
$\min\{n,k\kern-1pt+\kern-1pt1\}.$}
\vskip0pt}
\vskip-27pt\noi
}
}                                                       \eqno\hbox{\nrm(6.2)}
$$

\vskip20pt\noi
{\bf Remark 6.2.}\
A word on the separability assumption is in order.\ 
\vskip6pt\noi
(a)
It was applied in the original proof of (6.2) in \cite[Theorems 3.2]{SK}
that~the Hilbert space $\H$ upon which the normal operator $T^n\kern-1pt$ acts
is a direct integral$.$~Direct integral is a generalisation of direct sum of
Hilbert spaces, denoted by ${\H\kern-1pt=\!\int_X^{_\oplus}\H_x\kern1pt d\mu}$,
where $\{\H_x\}_{x\in X}$ is a field of Hilbert spaces and $\mu$ is any Borel
measure on a Borel $\sigma$-algebra $\X$ of subsets of a Borel index set
${X\!\sse\!\CC}.$

\vskip6pt
$\!$The common definition of direct integrals requires that the field
$\{\H_x\}_{x\in X}\!$ be~made up of separable Hilbert spaces (also referred to
by saying that the fibre spaces are separable).\ Nevertheless, it is also
common to assume that $\H$ is itself separable, as was the case in (6.2).\ For
a gentle summary on direct integrals, see, e.g., \cite[Section 4.1]{Geo} (see
also \cite{AC}$\kern.5pt$); for a full treatment, see
\cite[Part II, Chapters 1 and 2]{Dix}$.$ As introduced in
\cite[p$.\,$164,~Definition 1(iii), Section 3, Chapter$\;$1, Part$\;$II]{Dix},
the notion of direct integrals imposes~separability for the Hilbert spaces
$\H_x.$

\vskip6pt
The assumption that $\H$ is separable makes every $\H_x$ separable, as
required$.$~$\kern-1pt$The requirement that each $\H_x$ is separable, however,
does not imply that $\H$ is separa\-ble.\ But excludes Hilbert spaces $\H$
such as ${L^2(X,H,\mu)}$, the collection of all square-integrable functions
(with respect to a Borel measure $\mu$) of a Borel subset $X\kern-1pt$ of
$\CC$ into an arbitrary nonseparable Hilbert space $H.$ Reformulations for
fields of nonseparable Hilbert spaces have been investigated.\ For instance,
an approach towards a separability-free definition of direct integrals has
been carefully discussed~in~\cite{Wil}.\

\vskip6pt\noi
(b)
For the particular case where $\H$ is assumed to be separable {\it a priori}\/,
the direct integral of a normal operator (such as $\kern-.5pt T^n\kern-1pt$ in
(6.2)) with respect to a Borel measure $\mu$ is naturally identified via the
Spectral Theorem with the representation of the normal operator as a
multiplication operator on $L^2(X,\mu)$, where such~a~measure $\mu$ is finite
whenever $\H$ is separable and coincides with the scalar spectral measure of
the normal operator if it acts on a separable Hilbert space (see, e.g.,
\cite[Corollary 3.14, Lemma 4.7 and its Remarks]{ST2}$\kern.5pt).$

\vskip6pt
In fact, consider the version of the Spectral Theorem that says, ``if $N$ is a
normal operator on an {\it arbitrary}\/ complex Hilbert space $\H$, then there
is a measure $\mu$ on a $\sigma$-algebra of Borel subsets of the disjoint
union $\Omega$ of a collection $\{\Omega_\gamma\}$ of compact sub\-sets
$\Omega_\gamma$ of $\CC$ so that the identity function with
\hbox{multiplicity},~${\vphi\!:\Omega\to\CC}$~in~${L^\infty(\Omega,\mu)}$,
is such that $N$ is unitarily equivalent to the multiplication operator
$M_\vphi\kern-.5pt$ on $L^2(\Omega,\mu).\!$'' Moreover,
${\vphi(\Omega)^-\!=\sigma(N)}$, the spectrum of $N.$ $\kern-1pt$The proof of
it shows that, since $N$ on $\H$ is normal, $\H$ is the orthogonal direct sum
${\bigoplus_\gamma\!\M_\gamma}\kern-.5pt$ of subspaces $\M_\gamma\kern-.5pt$
of~$\H$,~where
$$
\M_\gamma={\bigvee}_{m,n}\!\big\{N^n N^{*m}x_\gamma\big\}
$$
for a family $\{x_\gamma\}$ of unit vectors in $\H.$ $\kern-1pt$Hence each
$\M_\gamma\kern-.5pt$ is spanned by a countable set so that
$\M_\gamma\kern-.5pt$ is separable (see, e.g., \cite[part (b), proof of
Theorem 3.11]{ST2}$\kern.5pt$).\ $\kern-1pt$Thus~if $N$ is normal on $\H$,
then ${\H\kern-1pt=\bigoplus_\gamma\!\M_\gamma}$, a (possibly uncountable)
direct sum of separable Hilbert spaces, which suggests that if $N$ is normal
on $\H$, then $\H$ is a direct integral.\ $\kern-1pt$There are measurability
conditions to be verified beyond the existence~of~a field $\{\M_\gamma\}$ of
separable Hilbert spaces.\ $\kern-1pt$The above measure $\mu$ equipping
$L^2(\Omega,\mu)$ is not necessarily finite.\ A sufficient condition for $\mu$
to be finite is the separability of $\H.$ $\kern-1pt$This separable case has
been considered in \cite[Section 3]{AC}, as we will see below.\

\vskip6pt\noi
(c)
Let $\kern-.5pt N\kern-.5pt$ be a normal operator on separable Hilbert space
$\H.$ So we may~\hbox{assume}~that
$\kern-.5pt\H\kern-1pt=\!{\int_X^{_\oplus}\!\H_\lambda d\mu}$ is a direct
integral.\ In this case, $N$ is given by
$\kern-.5pt N\!=\!{\int_X^{_\oplus}\!\lambda\,I_\lambda d\mu}$,
where~$\lambda$ is in ${ X\!=\sigma(N)}$, the spectrum of $N$, and each
$I_\lambda$ is the identity on $\H_\lambda.$ Recall~that~$\mu$ is a Borel
measure on the $\sigma$-algebra $\X$ of Borel subsets of the compact set
${X\kern-1pt=\sigma(N)}$, thus $\mu$ is a finite measure on~$\X.$ Since $\H$
is separable, the measure $\mu$ in item (b)~above is finite as well.\ (And
finite measures can be normalised into a probability measure.)\ If\/ $T^n$
is normal (as in (6.2)$\kern.5pt$) and $\H$ is separable, then
$T^n\!=\!{\int_{\sigma(T^n)}^{_\oplus}\!\lambda\,I_\lambda\,d\mu}.$ A crucial
property of the von Neumann algebra generated by $T^n$ is that every operator
in its commutant is represented as a decomposable operator on
$\H=\!{\int_{\sigma(T^n)}^{_\oplus}\!\H_\lambda\,d\mu}.$ $\kern-1pt$Thus
$T\kern-1pt=\!{\int_{\sigma(T^n)}^{_\oplus}\!T_\lambda\,d\mu}$ with each
$T_\lambda$ acting on each $\H_\lambda$, where $T\kern-1pt$ on $\H$ is
normal~if~and only if every operator $T_\lambda$ on $\H_\lambda$ is normal,
and $T^n\!=\!{\int_{\sigma(T^n)}^{_\oplus}\!T_\lambda^n\,d\mu}.$ Hence,
${T_\lambda^n=\lambda\,I_\lambda}.$

\vskip6pt
The programme in the proof of (6.2) was to ensure that if $T^n\!$ is normal
and $T$ is $k$-quasiparanormal, then $T_\lambda$ is paranormal.\ Since
${T_\lambda^n=\lambda I_\lambda}$ is normal as we saw above, $T_\lambda$ must
be normal by (5.3), and so~is~$T.$                       \hfill$\blacksquare$

\vskip6pt
Actually, the separability requirement in (6.2) is due to the application of
direct integrals in its proof.\ In Theorem 6.3 below, we remove the need for
separability, allowing arbitrary Hilbert spaces, by avoiding the usage of
direct integrals.\

\vskip6pt\noi
{\bf Theorem 6.3.}\
{\it If\/ $T\kern-1pt$ is\/ $k$-quasiparanormal, and\/ $T^n\!$ is normal,
then\/ ${T\kern-1pt=N\oplus L}$, where\/ $N$ is normal and\/ $L$ is
nilpotent of index at most\/ $\min\{n,k\kern-1pt+\kern-1pt1\}.$}

\proof
Let $\kern-.5pt T\kern-1pt$ be an operator on a Hilbert space $\kern-.5pt\H$,
and let $k$ and $n$ be arbitrary positive integers.\ Suppose $T\kern-1pt$ is
$k$-quasiparanormal and $T^n\kern-1pt$ is normal.\ Recall that if
$T^n\kern-1pt$ is normal, then $T^{m\kern.5pt n}\kern-1pt$ is normal for every
positive integer $m.$ Thus, with no loss of generality, assume throughout part
(a) of this proof that $T^n\kern-1pt$ is normal~for~some~${n\ge k}.$

\vskip6pt\noi
(a)
As we have observed before, it is clear that~$\R(T^n)$ is
$T\kern-1pt$-invariant, and so is $\R(T^n)^-\!.$ Since $T\kern-1pt$ commutes
with the normal operator $T^n\!$, the Fuglede Theorem ensures that
$T\kern-1pt$ commutes with $T^{n*}\!$ (see, e.g.,
\cite[Corollary 3.19]{ST2}$\kern.5pt).$ $\kern-1pt$Thus
${T\kern1pt T^{n*}\!=T^{n*}T}$ so that ${T^nT^*\!=T^*T^n}\!.$ Hence $\R(T^n)$
is $T^*\!$-invariant, and so is $\R(T^n)^-\!.$ $\kern-1pt$Therefore,
$$
\hbox{$\R(T^n)^-$ reduces $\,T.$}
$$
\vskip-3pt\noi
Thus,
$$
T\kern-1pt=T|_{\R(T^n)^-\!}\oplus T|_{\R(T^n)^\perp}
\quad\hbox{on}\quad
\H=\R(T^n)^-\!\oplus\R(T^n)^\perp,
$$
so that
$$
T^n\!=\big(T|_{\R(T^n)^-\!}\big)^n\!
\oplus\big(T|_{\R(T^n)^\perp\!}\big)^n\!.
$$
Moreover,
$$
\hbox{$\big(T|_{\R(T^n)^-\!}\big)^n$ is a normal operator}
$$
\vskip3pt\noi
since $T^n$ is normal.\ Furthermore, $T\kern-1pt$ is $n$-quasiparanormal
because it is $k$-quasipara\-normal and we have taken ${n\ge k}.$ (Recall:\ if
$T\kern-1pt$ is $k$-quasiparanormal, then $T\kern-1pt$ is $m$-quasiparanormal
for every integer ${m\ge k}.$) Thus, by Proposition 6.1(a), we get~that
\goodbreak\vskip4pt\noi
$$
\hbox{$T|_{\R(T^n)^-}$ is paranormal.}
$$
\vskip2pt\noi
Then, since $\big(T|_{^\R(T^n)^-\!}\big)^n$ is normal, it follows by (5.3)
that
$$
\hbox{$T|_{\R(T^n)^-}$ is normal}.
$$
As for the direct summand
${T|_{\R(T^n)^\perp}\!\!:\R(T^n)^\perp\!\to\R(T^n)^\perp}\!$, it is
necessarily nilpotent.\ Indeed, recall that $\R(T^n)^\perp\!=\N(T^{n*})$, so
that ${T^{n*}v=0}$ for every $v$ in $\R(T^n)^\perp\kern-2.5pt$,~and so
${T^{n*}|_{\R(T^n)^\perp}\!=O}$, which implies
${(T^n|_{\R(T^n)^\perp})^*\!=O}$ because $\R(T^n)^\perp\!$ reduces
$T\kern-1pt.$~So
$$
(T|_{\R(T^n)^\perp}\!)^n=T^n|_{\R(T^n)^\perp}\!=O.
$$
\vskip-2pt\noi
Consequently,
$$
T\kern-1pt=N\oplus L,
$$
where ${N\kern-.5pt=T|_{\R(T^n)^-}\!}$ is a normal operator and
${L=T|_{\R(T^n)^\perp}\!}$ is nilpotent such that ${L^n\!=O}$, and so $L$ is
a nilpotent operator of index ${j\le n}$ whenever ${n\ge k}.$

\vskip6pt\noi
(b)
To verify that the nilpotence index $j$ of $L$ is such that
${j\le\min\{n,k\kern-1pt+\kern-1pt1\}}$, we split the proof into two parts.
\vskip6pt\noi
(b$_1$)
The direct summand $L$ is $k$-quasiparanormal because $T\kern-1pt$ is
$k$-quasiparanormal.\ So $L$ is $m$-quasiparanormal for every integer
${m\kern-1pt\ge\kern-1pt k}.$ If ${L^{m+2}\!=O}$, then ${L^{m+1}\!=O}$ by the
definition of an $m$-quasiparanormal operator, for every
${m\kern-1pt\ge\kern-1pt k}.$ $\kern-1pt$Thus we can infer that if
${L^\ell\!=O}$ for some ${\ell\kern-1pt\ge\kern-1pt k\kern-1pt+\kern-1pt1}$,
then ${L^{k+1}\!=O}$ and so ${L^\ell\!=O}$ for all
${\ell\kern-1pt\ge\kern-1pt k\kern-1pt+\kern-1pt1}.$ Since $L$ is indeed
nilpotent, the nilpotence index $j$ of $L$~is such that
${j\le k\kern-1pt+\kern-1pt1}.$

\vskip6pt\noi
(b$_2$)
Moreover, since ${T^n\!=N^n\kern-1pt\oplus L^n}\kern-1pt$ is normal,
$L^n\kern-1pt$ is normal.\ But $L^n\kern-1pt$ is nilpotent because $L$ is
nilpotent.\ Hence ${L^n\kern-1pt=O}$ because the only nilpotent normaloid
operator is the null operator. So the nilpotence index $j$ of $L$ is such
that ${j\le n}.$                                                         \qed

\vskip6pt
It is clear that any of the direct summands in Theorem 6.3 may be absent~from
the decomposition of $T\kern-1pt.$ $\kern-1pt$Theorem 6.4 below extends
\cite[Theorem 3.5]{SK} to nonseparable Hilbert spaces.\

\vskip6pt\noi
{\bf Theorem 6.4.}\
{\it Let\/ $T\kern-1pt$ be a\/ $k$-quasiparanormal operator, and suppose\/
$T^2\kern-1pt$ is normal$.$ Consider the decomposition of\/ $T\kern-1pt$ in
Theorem 6.3.\ If\/ the nilpotent part\/ $L$ is nonzero, then it is given by\/
${L\!=\big(\smallmatrix{O & C \cr
                        O & O \cr}\big)}$
with\/ $C$ being an injective nonnegative operator}\/.

\proof
Let $T\kern-1pt$ be $k$-quasiparanormal and $T^2\kern-1pt$ normal.\ So
${T\kern-1pt=N\oplus L}$ by Theorem~6.3.\

\vskip6pt\noi
First recall that, since $L$ is nilpotent, it is pure (i.e., it has no normal
direct summand) because the~only normaloid nilpotent operator is the null
operator.\ Next, consider the following result from \cite[Theorem 1]{RR}.\
\vskip6pt\noi
{\narrower\narrower
{\it An operator is the square root of a normal operator if and only if it is
of the form\/
${A\oplus\big(\smallmatrix{B & \kern4pt C  \cr
                           O & \kern-1pt-B \cr}\big)}$,
where\/ $A$ and\/ $B$ are normal operators and\/ $C$ is positive\/ $($i.e.,
nonnegative and injective\/$)$ and commutes with\/ $B$.
\vskip0pt}
}
\vskip6pt\noi
Since $T^2\kern-1pt$ is normal, $L^2\kern-1pt$ is normal.\ Since $L$ is pure,
the above italicised result says that
${L\kern-1pt=\kern-1pt\big(\smallmatrix{B & \kern4pt C\kern-.5pt  \cr
                                        O & \kern-1pt-B\kern-.5pt \cr}\big)}$
with $C$ being a positive operator commuting with $B$,
where~$B$~is~\hbox{normal}.\ By Theorem 6.3, ${L\ne O}$ is nilpotent of index
not greater than ${n=2}.$ Thus
$L^2\!={\big(\smallmatrix{B^2 & O            \cr
                          O   & B^2\kern-1pt \cr}\big)=O}$
so that ${B^2\!=O}.$ Then ${B=O}$ because it is normal~and~\hbox{nilpotent}.\
                                                                         \qed

\vskip6pt
Under the assumption of quasiposinormality, again along the line in Lemma 5.2,
the nilpotent direct summand in Theorem 6.3 is always absent.\

\vskip6pt\noi
{\bf Theorem 6.5.}\
{\it Suppose\/ $T\kern-1pt$ is a quasiposinormal operator.\ If\/
$T\kern-1pt$ is\/ $k$-quasi\-paranormal and\/ $T^n\!$ is normal, then\/
$T\kern-1pt$ is normal}\/.

\proof
According to Theorem 6.3, if $T\kern-1.5pt$ is $k$-quasiparanormal on an
arbitrary \hbox{Hilbert} space $\H$ and $T^n\!$ is normal for some positive
integer $n$, then
$$
T\kern-1pt=N\oplus L
\quad\;\hbox{on}\;\quad
\H=\M^\perp\oplus\M
$$
for some subspace $\M$ of $\H$ that reduces $T\kern-1pt$, where $N$ is normal
and $L$ is nilpotent of index less than or equal to
$\min\{n,k\kern-1pt+\kern-1pt1\}.$ Since $T^n\kern-1pt$ is normal if and only
if $T^{m\kern.5pt n}\kern-1pt$ is normal for every positive integer $m$, take
again (with no loss of generality) an $n$ large enough so that ${L^n\!=O}$,
equivalently, ${\N(L^n)=\M}.$ If $T\kern-1pt$ is quasiposinormal, then so is
$L$ because ${\R(N^*)=\R(N)}.$ Thus
$$
\N(L)=\N(L^2)
$$
by Proposition 4.1(e).\ As we saw in the proof of Theorem 5.6, this implies
that
$$
\N(L)=\N(L^n)=\M
\quad\;\hbox{for every}\;\quad
n\ge1.
$$
\vskip2pt\noi
So ${L=O}$ because it acts on $\M.$ Therefore, ${T\kern-1pt=N\oplus O}$ is
normal.\                                                                 \qed

\vskip6pt
As in Corollary 5.5, the next result is a particular case of Theorem 6.5 with
a simpler statement.\

\vskip6pt\noi
{\bf Corollary 6.6.}\
{\it Suppose\/ $T\kern-1pt$ is posinormal.\ If\/ $T\kern-1pt$ is\/
$k$-quasiparanormal and\/ $T^n\kern-1pt$ is normal, then\/ $T\kern-1pt$ is
normal}\/.

\vskip6pt
We close the paper with a result on dense-range $k$-quasiparanormal operators,
which refers back to Question 5.3, Theorem 5.4, and Theorem 6.5.\ Recall that
dense-range operators are coquasiposinormal, and invertible operators are
posinormal and coposinormal.\

\vskip6pt\noi
{\bf Theorem 6.7.}\
{\it Let\/ $T\kern-1pt$ be\/ a $k$-quasiparanormal operator for some\/ $k.$}
\vskip4pt\noi
{\rm(a)}
{\it Suppose\/ $T\kern-1pt$ is invertible.\ $\kern-1pt$Then}
\begin{description}
\item{$\kern-11pt$\rm(a$_1$)$\kern2pt$}
{\it $T\kern-1pt$ and\/ $T^{-1}\kern-1pt$ are paranormal\/
$($and so they are\/ $k$-quasiparanormal for every\/ $k)$}.\
\vskip4pt
\item{$\kern-11pt$\rm(a$_2$)$\kern1pt$}
{\it If one of\/ $T^n\kern-1pt$ or\/ $T^{-n}\kern-1pt$ is quasinormal for
some\/ $n$, then so is the other,                                          \\
\hbox{$\kern4pt$and}\/ $T\kern-1pt$ and\/ $T^{-1}\kern-1pt$ are normal}.\
\end{description}
{\it More generally}\/,
\vskip4pt\noi
{\rm(b)}
{\it if\/ $T\kern-1pt$ has dense range and\/ $T^n\kern-1pt$ is quasinormal,
then\/ $T\kern-1pt$ is normal}\/.\

\proof
Let $T\kern-1pt$ be a $k$-quasiparanormal operator on a Hilbert space $\H.$

\vskip6pt\noi
(a$_1$)
According to Proposition 6.1(b), if an operator is invertible, then it is
$k$-quasiparanormal for an arbitrary $k$ if and only if it is paranormal.\
But an invertible paranormal operator $T\kern-1pt$ has a paranormal inverse
$T^{-1}\!$ \cite[Corollary 7.1.9]{Ist}.\
\vskip6pt\noi
(a$_2$)
If $T\kern-1pt$ is invertible, then $T^n$ is invertible.\ If the invertible
$T^n\kern-1pt$ is quasinormal, then $T^n$ is trivially normal.\ But
$T\kern-1pt$ is paranormal according to item (a$_1$).\ So $T$ is normal by
(5.3).\ Since an invertible operator and its inverse are quasinormal (thus
normal) together, the same conclusion comes out for $T^{-1}\!.$
\vskip6pt\noi
(b)
If $T\kern-1pt$ is $k$-quasiparanormal and ${\R(T)^-\kern-2.5pt=\kern-1pt\H}$,
then $T\kern-1pt$ is paranormal by Proposition 6.1.\ $\kern-1pt$Then apply
Theorem 5.4 since dense-range operators are coquasiposinormal$.\!\!$     \qed

\vskip-0pt\noi
\bibliographystyle{amsplain}

\begin{thebibliography}{10}

\bibitem{And}
T. Ando,
{\it Operators with a norm condition}\/,
Acta Sci. Math. (Szeged) {\bf 33} (1972), 169--178.
\goodbreak

\bibitem{AC}
E.A. Azoff and K.F. Clancey,
{\it Spectral multiplicity for integrals of normal operators}\/,
J. Operator Theory {\bf 3} (1980), 213--235.
\goodbreak

\bibitem{BKT}
P.S. Bourdon, C.S. Kubrusly, and D. Thompson,
{\it Powers of posinormal Hilbert-space operators}\/,
(2022) available at
https://arxiv.org/abs/2203.01473
\goodbreak

\bibitem{BKLT}
P.S. Bourdon, C.S. Kubrusly, T. Le, and D. Thompson,
{\it Closed-range posinormal operators and their products}\/,
Linear Algebra Appl. {\bf 671} (2023), 38--58.
\goodbreak

\bibitem{Con}
J.B. Conway,
{\it The Theory of Subnormal Operators}\/,
Mathematical Surveys and Monographs, Vol. 36, Amer. Math. Soc.,
Providence, 1991.
\goodbreak

\bibitem{Dix}
J. Dixmier,
{\it Von Neumann Algebras}\/,
North-Holland, Amsterdam, 1981;
translation of {\it Les Alg\`ebres d'Op\'erateurs dans l'Espace Hilbertien
$($Alg\`ebres de Von Neumann\/$)$}\/,
2\,\`eme \'ed., Gauthier-Villars, Paris, 1969 (1\,\`ere \'ed., 1957).
\goodbreak

\bibitem{Dug}
B.P. Duggal,
{\it On $n$th roots of normal contractions}\/,
Bull. London Math. Soc. {\bf 25} (1993), 74--80.


\bibitem{Geo}
A. Geondea,
{\it When are products of normal operators normal}\/?
Bull. Math. Soc. Sci. Math. \hbox{Roumanie} {\bf 52}(100) (2009), 129--150.
\goodbreak

\bibitem{Ito}
M. Itoh,
{\it Characterization of posinormal operators}\/,
Nihonkai Math. J. {\bf 11} (2000), 97--101.
\goodbreak

\bibitem{Ist}
V.I. Istr\v a\c tescu,
{\it Introduction to Linear Operator Theory}\/,
Marcel Dekker, New York, 1981.
\goodbreak

\bibitem{JKKP}
I.H. Jeon, S.H. Kim, E. Ko, and J.E. Park,
{\it On positive-normal operators}\/,
Bull. Korean Math. Soc. {\bf 39} (2002), 33--41.
\goodbreak

\bibitem{Ker}
L. K\'erchy,
{\it On roots of normal operators}\/,
Acta Sci. Math. (Szeged) {\bf 60} (1995), 439--449.

\bibitem{HSO}
C.S. Kubrusly,
{\it Hilbert Space Operators:\ A Problem Solving Approach}\/,
Birkh\"auser, Boston, 2003.
\goodbreak

\bibitem{ST2}
C.S. Kubrusly,
{\it Spectral Theory of Bounded Linear Operators}\/,
Birkh\"auser-Springer-Switzer\-land, Cham, 2020.
\goodbreak

\bibitem{KD}
C.S. Kubrusly and B.P. Duggal,
{\it On posinormal operators}\/,
Adv. Math. Sci. Appl. {\bf 17} (2007), 131--148.
\goodbreak

\bibitem{KVZ}
C.S. Kubrusly, P.C.M. Vieira, and J. Zanni,
{\it Powers of posinormal operators}\/,
Oper. Matrices {\bf 10} (2016), 15--27;
{\it Erratum/Addendum}\/, Oper. Matrices {\bf 16} (2022), 1239--1242.
\goodbreak

\bibitem{PS1}
P. Pietrzycki and J. Stochel,
{\it Subnormal\/ $n$th roots of quasinormal operators are quasinormal}\/,
J. Funct. Anal. {\bf 280}(12) (2021), $\#$109001, 14 pp.
\goodbreak

\bibitem{PS2}
P. Pietrzycki and J. Stochel,
{\it On nth roots of bounded and unbounded quasinormal operators}\/,
Ann. Mat. Pura Appl. {\bf 202} (2023), 1313--1333.
\goodbreak

\bibitem{RR}
H. Radjavi and P. Rosenthal,
{\it On roots of normal operators}\/,
J. Math. Anal. Appl. {\bf 34} (1971), 653--664.
\goodbreak

\bibitem{Rha}
H.C. Rhaly, Jr.,
{\it Posinormal operators}\/,
J. Math. Soc. Japan {\bf 46} (1994), 587--605.
\goodbreak

\bibitem{JV}
P. Sam Johnson and A. Vinoth,
{\it Product and factorization of hypo-EP operators}\/,
Spec. Matrices {\bf 6} (2018), 376--382.
\goodbreak

\bibitem{Stampfli}
J.G. Stampfli,
{\it Hyponormal Operators}\/,
Pacific J. Math. {\bf 12} (1962), 1453--1458.
\goodbreak

\bibitem{Sta}
H. Stankovi\'c,
{\it Subnormal $n$-th root of matricially and spherically quasinormal pairs}\/,
Filomat {\bf 3}(16) (2023), 5325--5331.
\goodbreak

\bibitem{Sta2}
H. Stankovi\' c,
{\it Spherically quasinormal tuples:\ $n$-th root problem and hereditary
properties}\/,
Complex Anal. Oper. Theory, {\bf 19}(6) (2025), 156:1--16.

\bibitem{SK}
H. Stankovi\'c and C.S. Kubrusly,
{\it On roots of normal operators and extensions of Ando's Theorem}\/,
Ann. Funct. Anal. {\bf 16}(4) (2025), 60:1--15.
\goodbreak

\bibitem{Wil}
W. Wills,
{\it Direct integrals of Hilbert spaces I}\/,
Math. Scand. {\bf 26} (1970), 73--88.
\goodbreak

\bibitem{YY}
T. Yamazaki and M. Ynagida,
{\it Relations between two inequalities and their application to paranormal
operators}\/,
Acta Sci. Math. (Szeged) {\bf 69} (2003), 377--389.
\goodbreak

\end{thebibliography}

\end{document}